\documentclass{article}
\usepackage{amsmath,amsfonts,amsthm,amssymb,amscd,colordvi}
\newcommand{\R}{{\mathbb R}}
\newcommand{\C}{{\mathbb C}}
\newcommand{\Prr}{{\mathbb P}}
\newcommand{\p}{{\partial}}
\newcommand{\e}{\varepsilon}
\newcommand{\id}{\mathop{\rm id}\nolimits}

\newcommand{\vp}{\varphi}

\newcommand{\const}{\mathop{\rm const}\nolimits}

\newcommand{\Z}{{\mathbb Z}}

\newcommand{\T}{{\mathbb T}}
\newcommand{\cH}{{\cal H}}

\newcommand{\N}{{\mathbb N}}

\theoremstyle{plain}
\newtheorem{theorem}{Theorem}[section]

\newtheorem{corollary}[theorem]{Corollary}
\theoremstyle{definition}

\theoremstyle{remark}

\newtheorem{conjecture}[theorem]{Conjecture}

\begin{document}

\author{Sergei B. Kuksin\,$^1$ and Anatoly I. Neishtadt\,$^2$}
\footnotetext[1]{ CNRS and CMLS at Ecole Polytechnique, 
Palaiseau, France; 
 e-mail: {\tt kuksin@math.polytechnique.fr}}
 \footnotetext[2]{ Loughborough University, United Kingdom and Space Research Institute, Moscow, Russia;
  e-mail: {\tt  A.Neishtadt@lboro.ac.uk}}

\title{On quantum averaging, quantum {KAM} and quantum diffusion
}
\date{}
\maketitle

\hfill{\it Dedicated to the memory of Mark Iosifovich Vishik. } \medskip

\begin{abstract}
For nonautonomous Hamiltonian systems and their quantisations we discuss 
properties of the quantised systems, related to those of 
the corresponding classical systems,  described by the KAM-related theories: 
 the proper KAM, the averaging theory, the Nekhoroshev stability, and 
the diffusion.
\end{abstract}

\section{Introduction}\label{s:intr}
Consider a classical nonautonomous  Hamiltonian system on the phase-space $T^*\T^d=\R^d\times\T^d=\{(p,q)\}$ 
or  $T^*\R^d=\R^d\times\R^d$ 
with a Hamiltonian $H(p,q,t)$:
\begin{equation}\label{0.1}
 \dot p=-\nabla_q H,\qquad \dot q=\nabla_p H.
\end{equation}
The corresponding quantum Hamiltonian operator is obtained by replacing in $H(p,q,t)$ the variable
$q_j$, $j=1,\dots,d$, by the operator which acts on complex functions $u(x)$ as 
multiplying by $x_j$, and replacing each $p_j$ by the operator $\frac{\hbar}{i} \frac{\p}{\p x_j} $, where 
$\hbar$ is the Planck constant.\footnote[3]{This rule  of quantisation is the most common, but certainly it is not
unique. More generally one  may replace $q_j$ and $p_j$ by any operators $Q_j$  and $P_j$ such that
$[Q_j,P_k]=i\hbar\delta_{j,k}$, for all $j$ and $k$.}
The Hamiltonian operator 
$
\cH=H(\frac{\hbar}{i} \nabla_x,x,t)
$
defines a quantum system, and a classical problem of the quantum mechanics, streaming from its
first years of existence, is to study (spectral) properties of the operator $\cH$ and the properties of the 
corresponding evolutionary equation 
\begin{equation}\label{0.2}
 i\hbar\, \dot u(t,x)=  \cH u(t,x),
\end{equation}
in their relation with the classical  system \eqref{0.1}.

For example, if 
\begin{equation}\label{0.3}
 H(p,q,t)=|p|^2+V(t,q),
\end{equation}
then
\begin{equation}\label{0.0}
\cH=\cH_t=-\hbar^2\Delta+V(t,x),
\end{equation}
i.e. $\cH$ is the Schr\"odinger operator with the potential $V$. 

In this paper we discuss properties of the Hamiltonian operator $\cH$, corresponding to properties of 
system \eqref{0.1},  described by the KAM-related theories. Namely, by 
 the proper KAM, the averaging,  the Nekhoroshev stability, and 
the diffusion (this list by no means is canonical; it corresponds to the authors'  taste).  We 
 discuss results for quantum systems \eqref{0.2}  which we regard as parallel to the three classical 
theories above, mostly restricting ourselves to the case of periodic boundary conditions $x\in\T^d$ and
 assuming that $\hbar=\,$const. Scaling $x$ and $t$ in the dynamical equation \eqref{0.2},  \eqref{0.0} we
achieve $\hbar=1$. A discussion concerning semiclassical limit  $\hbar\to0$, when  it is not appropriate
to scale $\hbar$ to 1,  is contained in Section~\ref{quasi-classic}. There we consider the equations in the whole
space, $x\in\R^d$, since for the periodic boundary conditions the corresponding results are less developed.

All quantum results we discuss deal with non-autonomous equations \eqref{0.2}, \eqref{0.0},
  so  their classical analogies  are ``KAM-related'' theories for non-autonomous Hamiltonian systems \eqref{0.3}. 
We do not touch very interesting, important and complicated problem of constructing eigenfunctions of
nearly integrable Hamiltonian operators by quantasing KAM-tori of the corresponding  autonomous
Hamiltonian  systems 
(see \cite{Laz}). 

\medskip

Let $u(t)$ be a solution of the equation \eqref{0.2}, \eqref{0.0}. Multiplying the equation by $\bar u$ and integrating
 over $\T^d$ we get that
 $\ 
    |u(t)|_{L_2}^2=\const.
 $
 Write $\ u(t,x)=\sum_su_s(t) \vp_s(x)$, where $\{\vp_s\}$ are eigenfunctions of the ``unperturbed" Hamiltonian operator.
 Then  $\sum |u_s(t)|^2\equiv\,$const. What happens to the quantities 
 $|u_s(t)|^2$ as $t$ growths, i.e. how the total probability $\sum |u_s(t)|^2$  is distributed between the states $s\in\Z^d$
 when $t$ is large? This is the question which is  addressed by the theorems we discuss.
 \medskip

 \noindent
 {\bf Acknowledgments.} The authors are thankful for discussions to  Sergey Dobrokhotov,  H{\aa}kan~Eliasson, and  Johannes Sj\"ostrand.  SK acknowledges the support  of  l'Agence Nationale de la Recherche through the grant
  ANR-10-BLAN~0102.

 \section{Quantum averaging}
 \subsection{Averaging and adiabatic invariance}\label{ss1}
 Let a classical Hamiltonian \eqref{0.3} have the form 
 \begin{equation}\label{ad_class}
 H(p,q,\e t)=H_{\e}=|p|^2+  V(\e t,q), 
\end{equation}
where the unperturbed Hamiltonian $|p|^2+  V(\tau,q)\,,\ \tau={\rm const}$, is integrable for each 
 $\tau$. Let $I_j, 1\le j\le d$, be the corresponding actions.
The classical  averaging principle  (e.g., see in \cite{AKN, LochM}) implies that  each action  is an  adiabatic invariant, namely   if $u_\e(t)$ is a solution of 
the perturbed equation \eqref{0.1}${}_{H=H_\e}$,
then $I_j(u_\e(t))$ stays almost constant on time-intervals of order $\e^{-1}$ . The averaging principle  is a  heuristic statement, and it does not always lead to correct results.  The adiabatic invariance for classical systems is discussed in more details in Section~\ref{quasi-classic}.

Now let us drop the assumption that the Hamiltonians \eqref{ad_class} with frozen $t$ are integrable and 
 consider the   corresponding 
  quantum system: 
\begin{equation}\label{S}
 \dot u=   -i\big(-\Delta u+ V (\e t,x) u\big),\quad 
  x\in\T^d.
\end{equation}
We assume that the function   $V(\tau,x)$ is  $C^2$-smooth bounded and 
denote by $A_{\e t}$ the linear operator in \eqref{S}, 
$$A_{\e t}=-\Delta + V(\e t,x). 
$$
Let  $\{\vp_s(\tau), s\in\Z^d\}$ and $\{\lambda_s(\tau)\}$ be the eigenvectors and the eigenvalues of $A_\tau$, where 
each $\lambda_s(\tau)$ is continuous in $\tau$. 
Let $u(t,x)$ be a solution of \eqref{S}, equal at $t=0$ to a pure state,
\begin{equation}\label{15}
u(0,x)= \vp_{s_0}(0),
\end{equation}
 such that for each $\e t$, $\lambda_{s_0}(\e t)$ is an isolated eigenvalue of $A_{\e t}$ of a constant multiplicity. Consider expansion of  $u(t,x)$ over the basis  $\{\vp_s(\tau), s\in\Z^d\}$:
$$
u(t,x)=\sum_su_s(t)\vp_s(\e t)\, .
$$ 
The quantum adiabatic theorem says that  $u(t,x)$ stays close  to the eigenspace, corresponding to $\lambda_{s_0}(\e t)$:

\begin{theorem}\label{tB-F}   (M.~Born, V.~Fock \cite{BF28} and T.~Kato \cite{Kat50}
) 
\begin{equation}\label{kato}
\sup_{0\le t\le \e^{-1}   }\sum_{s:\,  \lambda_s(\e t)\ne\lambda_{s_0}(\e t)}|u_s(t)|^2\to0\;\;\text{as}\;\; \e\to0. 
\end{equation}
\end{theorem}

This is a very general result which remains true for systems in the whole space (when $x\in\R^d$) 
 if  the operators
$A_{\e t}$ have  mixed spectrum, 
but $\lambda_{s_0}(\e t)$ always is an isolated eigenvalue of constant multiplicity, see in \cite{LochM}.
The case when this eigenvalues may be approached by other eigenvalues is considered in \cite{AvEl}.

Both for  classical and quantum systems,     adiabatic theorems   are considered often   on infinite time interval $-\infty<t<\infty$ under condition that the dependence of the potential $V$ on time disappears fast enough  as   $t\to\pm\infty$, and the system is sufficiently smooth.   In this case, for classical   Hamiltonians with $d=1$,     difference of values of actions on  trajectory at $t\to\pm\infty$ tends to 0 much faster than $\e$ as $\e\to 0$;   in  analytic case this difference is  $O(\exp(-{\rm const}/\e))$,  see \cite{LL} and references in  \cite{AKN},  Sect. 6.4.5.  For quantum systems, if  for $\tau\to-\infty$
all the probability is concentrated 
in the states, corresponding to the eigenvalue $\lambda_{s_0}(\tau)$,   then all the probability but a very small remnant will be absorbed by these states   as $\tau\to+\infty$. In  analytic case this remnant  is  $O(\exp(-{\rm const}/\e))$ \cite{Ne, JP} (this result also follows  from the calculus, developed 
in \cite{S}). 

We will return to the quantum adiabaticity  in Section~\ref{quasi-classic}. We note that also 
 there are  adiabatic theorems for systems where the Hamiltonian slowly depends not only on time, but also 
 on a part of the space-variables, e.g. see \cite{AKN}, Sect. 6.4.1 for classical systems and  \cite{Dobr} for quantum systems.

 \subsection{Around  Nekhoroshev's Theorem.
 }\label{ss2}
 Let us 
 start with classical systems. 
   Let 
 $H_\e(p,q)=h_0(p)+\e h_1(p,q)$, where the function 
  $h_0$ is analytic and steep (e.g., it is strictly 
 convex, for the definition of steep functions see \cite{Nek77} and \cite{LochM, AKN}). 
 Let $(p(t),q(t))$ be a solution of \eqref{0.1}. 
 Then there are $a,b>0$
 such that 
\begin{equation}\label{Nek}
|p(t)-p(0)|\le\e^a \qquad \forall\,|t|\le  e^{\e^{-b}},
\end{equation}
see in  \cite{Nek77, LochM, AKN}. 
There are many related results.  For example: let 
 $$
 H_\e(p,q,t)=h_0(p)+\e h_1(\omega t; p,q),\quad \omega\in\R^N,
 $$
where  $h_1$  is an  analytic function on $\T^N\times\R^d\times\T^d$,  $N\ge1$.
Then for a typical $\omega$ estimate \eqref{Nek} is true. 
 In particular, let us  take
 $$
 H_\e(p,q,t)=|p|^2+\e V(\omega t;q).
 $$
 The corresponding quantised Hamiltonian   is the operator
 $
 -\Delta +\e V(\omega t; x),
 $
 and the  evolutionary equation is 
 \begin{equation}\label{0}
 \dot u= - i \big(-\Delta u+\e V(\omega t;x) u\big).
 \end{equation}
 Do we have for solutions of \eqref{0}   an analogy of the Nekhoroshev estimate \eqref{Nek}?  I.e. is it true that 
 actions of the unperturbed system, evaluated along solutions of the perturbed equation \eqref{0}, do not change 
 much during exponentially long time? It turns out that a weaker form of this assertion holds true, even when $\e=1$! 
 Let us  consider the equation
\begin{equation}\label{0.02}
\dot u= -i \big(-\Delta u+ V(t, x) u\big),
\end{equation}
and 
consider the squared $r$-th Sobolev norm of $u$:
 $$
 \|u\|_r^2 = \sum_{s\in\Z^d} |u_s|^2(1+|s|^2)^r, \qquad r\in\R.
 $$
This is a linear combination of the actions for the unperturbed system with $V=0$.

\begin{theorem}\label{tB1}  (\cite{Bo99}). 
Let  
$V(t,x)=\tilde V(\omega t,x)$,  where $\omega\in\R^N$ is a Diophantine vector and  $\tilde V$ is a smooth 
function on $\T^N\times \T^d$. Then for each $r\ge1$ there exists $c(r)$
such that  any solution $u(t)$ of \eqref{0.02} satisfies
\begin{equation} \label{0.11}
\|u(t)\|_r\le \rm{Const}\cdot
 (\ln  t)^{c(r)}\|u_0\|_r,\qquad \forall\,t\ge2. 
\end{equation}
\end{theorem}

So if $u_0$ is smooth, then the high states $u_s$ stay almost non-excited for very long time. We miss a result
which would imply that the quantity in \eqref{kato}, calculated for solutions of \eqref{0.02}, \eqref{15} stays small for long time.

It is surprising that a weaker version of this result holds for potentials $V$ which are not time-quasiperiodic:
\begin{theorem}\label{tB2} (\cite{Bo99c}).
Let $V$ be smooth and  $C^k$-bounded uniformly in $(t,x)$ for each $k$. Then for each $r\ge1$ and $a>0$ 
there exists $C_a$
such that 
\begin{equation*} 
\|u(t)\|_r\le C_a\,    t^a   \|u_0\|_s,
\qquad \forall\,t\ge2. 
\end{equation*}
\end{theorem}

Also see \cite{Del10}. 
If  the potential $V(t,x)$ is analytic, then the norm $\|u(t)\|_r$ satisfies \eqref{0.11}, 
see \cite{WM08}. 
We are not aware   of any classical analogy of these results.

 \section{
 Quantum KAM }\label{s2}
 Let $(p,q)\in \R^d\times \T^d$.  Consider integrable Hamiltonian $h_0(p)=|p|^2$ and its time-quasiperiodic
  perturbation 
 $H_\e(p,q)=h_0(p)+\e V(\omega t,q)$, $\omega\in \R^n$, where 
 $V$ is analytic.  For the corresponding Hamiltonian
 equation we  have a KAM  result:
 {\it For a typical $(p(0),q(0))$ and a typical $\omega$ the solution $(p(t),q(t))$ is time-quasiperiodic. 
  }

 The quantised Hamiltonian defines the dynamical equation \eqref{0}. 
 We regard the vector $\omega$ as  a parameter of the problem:
$
\omega\in U\Subset \R^n.
$
We abbreviate  $L^2=L^2(\T^d,\C)$ and provide this space with the exponential basis
 $$\{e^{is\cdot x},s\in\Z^d \}.
 $$
For any linear operator  $B:L^2\to L^2$  let  $(B_{ab}, a,b\in\Z^d)$
 be its matrix in this basis. 
 
 The theorem below may be regarded as a quantum analogy of the KAM theorem above.
 For $d=1$ it is proven in  \cite{BG01}, and for $n\ge2$ --  in  \cite{EK09}. We do not know how to pass in this result to 
 the semiclassical limit.

\begin{theorem}\label{tEK2}
 If $\e\ll1$, then for most $\omega$ we can find an $\vp$-dependent
complex-linear isomorphism $\Psi(\vp)=\Psi_{\e,\omega}(\vp)$, \ $\vp\in\T^N$,
$$
\Psi(\vp):L^2\to L^2,\quad  u(x)\mapsto \Psi(\vp)u(x),
$$
and a bounded Hermitian operator $Q=Q^{\e,\omega}$ such that a curve $u(t)\in L^2$ solves eq.~\eqref{0}
if and only if $v(t)=\Psi( t\omega)u(t)$ satisfies 
$$
\dot v= i\big( \Delta v- \e Qv\big).
$$
The matrix $(Q_{ab})$ is block-diagonal, i.e. 
$\ 
Q_{ab}=0\quad\text{if}\quad |a|\ne|b|
$, and it satisfies
$$
Q_{ab}=(2\pi)^{-n-d}\int_{\T^N}\int_{\T^d} V(\vp,x)e^{i(a-b)\cdot x  }\,dxd\vp+O(\e^\gamma),\quad \gamma>0.
$$
Moreover, for any $p\in\N$ we have
$
\|Q\|_{H^p,H^p}\le C_1$ and   $ \|\Psi(\vp)-\id  \|_{H^p,H^p}\le \e C_2
$.
\end{theorem}

Here
``for most" means ``for all $\omega\in U_\e\subset U$, where mes$\,(U\setminus U_\e)\le\e^\kappa$
for some $\kappa>0$".  In particular, for any $\omega$ as in  the theorem all solutions of eq.~\eqref{0} are almost-periodic
functions of time. 
 Their Sobolev norms are almost constant:

\begin{corollary} For  $\omega$ as in the theorem and for
any $p$ solutions of \eqref{0} satisfy
$$
(1-C\e)\|u(0)\|_p\le \|u(t)\|_p\le (1+C\e)\|u(0)\|_p,\quad \forall\,t\ge0.
$$
\end{corollary}
This property  is  called the {\it dynamical localisaton}.

{\bf Proof.} Since $Q$ is block-diagonal, then $\|v(t)\|_p=\,$const. Since $v(t)=\Psi(t)u(t)$
and $\|\Psi-\id\|_{H^p,H^p}\le\e C_2$, then the estimate follows. \hfill $\Box$
\medskip

{\bf
Remarks.} 1) Let $n=0$. Then \eqref{0} becomes the equation 
$\
\dot u=-i \big(\Delta u+\e V(x) u\big).
$
Theorem states that this equation may be reduced to a block-diagonal equation
$\ \ 
\dot u=-i Au$, where\ $A_{ab}=0\;\;\text{if}\;\; |a|\ne|b|.\ 
$
This is a well known fact.

2) For $n=1$ the theorem's assertion is the Floquet theorem for the time-periodic equation \eqref{0}.
In difference with the finite-dimensional case, this is a perturbative result, valid only for  `typical' frequencies $\omega\in\R$ and small $\e$.

{\it 
 Proof of the Theorem}. Eq. \eqref{0} is a non-autonomous linear Hamiltonian system in $L^2$:
$$
\dot u=-i \frac{\delta}{\delta \bar u} H_\e(u),\quad H_\e(u)=\frac12 \langle\nabla u,\nabla \bar u\rangle
+\frac12 \e  \langle   V(\vp_0+t\omega,x)u,\bar u  \rangle.
$$
Consider the extended phase-space $L^2\times \T^n\times \R^n=\{(u,\vp,r)\}$.  There
 the  equation above can be written as the autonomous Hamiltonian  system
\begin{equation*}\begin{split}
&\dot u=-i \frac{\delta}{\delta \bar u} h_\e(u,\vp,r),\\
&\dot\vp=\nabla_r h_\e=\omega,\\
&\dot r=-\nabla_\vp h_\e,
\end{split}
\end{equation*}
where 
$\
h_\e(u,\vp,r,\e)=\omega\cdot r+
\frac12 \langle\nabla u,\nabla \bar u\rangle
+\frac12 \e  \langle   V(\vp,x)u,\bar u  \rangle.
$
So 
$h_\e$ is a small perturbation of the integrable quadratic Hamiltonian 
$\ h_0=\omega\cdot r+
\frac12 \langle\nabla u,\nabla \bar u\rangle
$.
To perturbations of $h_0$ applies the KAM-theorem from \cite{EK10}.
To show how this  implies the  Theorem~\ref{tEK2} let us write $h_\e$ as
$$
h_\e(u,\vp,r,\e)=\omega\cdot r+
\frac12 \langle\nabla u,\nabla \bar u\rangle
+\e f(u,\vp,r).
$$
In our case
$f=\frac12   \langle   V(\vp,x)u,\bar u  \rangle$. The theorem below is the main result of \cite{EK10}).

\begin{theorem}\label{tEK1}
 There exist a domain 
$\ \ {\cal O} =\{\|u\|<\delta\} \times \T^n\times\{|r|<\delta\}
$
and a symplectic transformation   $\ \ \Phi:{\cal O}\to L^2\times \T^n\times \R^n\ $
which transforms  $h_\e$ to
$$
h_0=\omega'\cdot r+\frac12\langle \nabla u,\nabla \bar u\rangle +\e \langle Qu,\bar u\rangle
+f'(u,\vp,r),
$$
where $f'=O(|u|^3)+O(|r|^2)$.
\end{theorem}

Torus$\ T_0=0\times \T^n\times 0\ $ is invariant for the transformed system, so
$\Phi(T_0)$ is invariant for the original equation. This is the usual KAM statement.
 Now it is trivial since it simply states that $u(t)\equiv 0$ is a solution on the original equation.

But the KAM theorem above tells more. Simple analysis of the proof (see a Remark in [EK2]) shows that if the
 perturbation $\e f$ is quadratic in $u$ and $r$-independent, then
 the KAM-transformations are linear in $u$ and do not change $\omega$.
 So the transformed Hamiltonians stay quadratic in $u$. Hence, the 
  Hamiltonian $h_0$ is such that  $f'=0$. That is, \\
$$
h_0=\omega'\cdot r+\frac12\langle \nabla u,\nabla \bar u\rangle +\e \langle Qu,\bar u\rangle.
$$
  This proves Theorem \ref{tEK2}.

\section{Quantum diffusion.}\label{s3}

 Let $(p,q)\in R^d\times \T^d$.  Consider 
 $H_\e(p,q)=|p|^2+\e V(\omega t,q)$,  where $\omega\in \R^N$ and $V$ is analytic.  Then

i) by KAM, for a typical $\omega$ and  typical initial data $(p_0,q_0)$ the solution such that $(p(0), q(0))=(p_0, q_0)$ is time-quasiperiodic;

ii) for exceptional $\omega$ and $(p_0,q_0)$ we ``should" have the Arnold diffusion: the 
action $p(t)$ of a corresponding solution  slowly 
``diffuses away"
from $p_0$.
 
   As before, 
 the quantised Hamiltonian defines the dynamical equation \eqref{0}.
 
{\bf Claim 4.1.}
Let $d=1$,  $N\ge2$ and the potential $V$ is nondegenerate in a suitable sense. 
Then there exist a smooth function $u(0,x)$  and $\omega\in \R^N$
such that 
\begin{equation}\label{growth}
\limsup_{t\to\infty} \| u(t)\|_s=\infty 
\end{equation}
for some $s\ge1$. 
\smallskip

An {\it example} of a time-periodic potential $V$, satisfying \eqref{growth}, is given in \cite{Bo99}. It is conjectured by
H.~Eliasson that the validity  of the Claim for a {\it typical} potential follows from the method of his work \cite{E02}.  
Proof of this assertion  is a work under preparation.

\section{Perturbed harmonic and anharmonic oscillators.}
In  Sections~\ref{s2},~\ref{s3} we
 deal with the evolutionary Schr\"odinger equation under periodic boundary 
conditions. Some similar results are available for equations in the whole space with  growing 
potentials:
\begin{itemize}

 \item  Consider Schr\"odinger equation in $\R^1$:
$$
\dot u= -i\big(-u_{xx} +(x^2+\mu x^{2m})u +\e V(t\omega,x)u\big),
$$
where $\mu>0, \ m\in\N,\ m\ge2$; $V(\vp,x)$ is $C^2$-smooth in $\vp, x$  and analytic in $\vp$, bounded uniformly 
in $\vp,x$.
An analogy of Theorem~\ref{tEK2} holds. See  \cite{K1} (Section 2.5) for the needed KAM-theorem.

 \item Due to Bambusi-Graffi \cite{BG01},   the  result holds for non-integer $m$. That is, for equations
$$
\dot u=-i \big( -u_{xx} +Q(x)u +\e V(t\omega,x)u\big),
$$
where $Q(x)\sim |x|^\alpha, \alpha>2$ as $|x|\to\infty$.  The potential $V$ may grow to infinity as  $|x|\to\infty$.

\item  Liu-Yuan  \cite{LY10}  allow faster growth of $V(x)$ in $x$.
Their result  applies to prove an  analogy of Theorem~\ref{tEK2} for  the {\it quantum Duffing oscillator }
$$
\dot u= - i\big(
-u_{xx} +x^4u +\e xV(t\omega,x)u\big).
$$
\item Due to Grebert and Thomann \cite{GT11},  the assertion holds for the  perturbed harmonic oscillator
$$
\dot u=-
i\big(-u_{xx} +x^2 u +\e V(t\omega,x)u\big).
$$
\end{itemize}

What happens in higher dimensions,  $d\ge2$ ? -- This is completely unknown. 

\section{Quantum adiabatic theorem in semiclassical limit}
\label {quasi-classic}
In this Section  we consider the classical system on $T^*\R^d=\R^d\times\R^d$ 
with a Hamiltonian
 \begin{equation}\label{ad_class1}
 H(p,q,\tau)=|p|^2+  V(\tau,q),
 \quad \tau =\e t,
\end{equation}  
and the corresponding quantum system
\begin{equation}\label{S_h}
  i\hbar \, \dot u=   -{\hbar}^2\Delta u+ V (\tau,x) u = \cH_\tau u,
  \quad \tau =\e t, 
\end{equation}
 (see \eqref{0.0}). We assume that for each $\tau$ the potential 
  $V(\tau,x)$ grows to infinity with 
$|x|$, so  the operator $\cH_\tau$ has  a discrete spectrum. 

We  fix small enough $\e$ that allows to make some statements about the
dynamics of the classical system, and then pass to the limit as $\hbar\to 0$. This limiting
 dynamics  may be  quite different from that in Section~\ref{ss1}
   when $\hbar$ is fixed and $\e\to 0$, as it was demonstrated  by M. Berry \cite{Ber84}  in the following striking example. Let $d=1$ and potential $V$ for $\tau={\rm const}$ has two (non-symmetric) potential wells.   Generically, for  $\tau={\rm const}$  and small enough     $\hbar$    each well supports a family of pure  quantum states localised mainly in this well.   Consider a
   solution $u(t,x)$ of equation (\ref{S_h}) with an initial condition which is a pure  quantum state from the left well.  For however  small $\e$  there exists $\hbar_0=\hbar_0(\e)>0$ such that if $0<\hbar<\hbar_0$, then for each $t\in[0,1/\e]$ the function  $u(t,\cdot)$  is localised in the same left well.  On the other hand, under some rather general assumptions, for however  small  $\hbar$   there exist $\e_0=\e_0(\hbar)$ and positive constants $ a_1< a_2$,  such that if $0<\e<\e_0$ then  the function  $u(t,\cdot)$ is localised in the right well for $  a_1\hbar/\e\le t\le a_2\hbar/\e$.
   
 Discussion of the case $\e\sim \hbar$ is contained in \cite{Kar90}. In what follows $\e_0, c, c_i$ are positive constants.

\subsection{Systems with one degree of freedom}
Assume first that 
 classical Hamiltonian  (\ref{ad_class1}) has one degree of freedom. We suppose that  $V$ is 
 $C^{\infty}$-smooth   and that in the phase plane of the Hamiltonian system  (\ref{ad_class1}) for each   $\tau={\rm const }$ there is a domain  filled  by closed trajectories. In this domain we  introduce action-angle variables  $I=I(p,q,\tau), \chi=\chi(p,q,\tau)\  {\rm mod}\ 2\pi$  (i.e.  $\chi\in\T^1$).    Invert these relations: $p=p(I, \chi, \tau), q=q(I, \chi, \tau)$. 
  Suppose that there is an interval $[a_1, b_1]$,   $ 0<a_1< b_1$, such that
  the map $(I, \chi, \tau) \mapsto (p, q, \tau)$ is smooth for  $I\in [a_1, b_1],\chi\in \T^1, \tau\in[0,1] $.
  We  express Hamiltonian (\ref{ad_class1}) via the 
  action variable and slow time: $H(p,q,\tau)=E(I, \tau)$. 

For  $\e >0$  let $(p(t),  q(t))$ be a solution of the perturbed  system with the Hamiltonian $H(p,q,\e t)$. 
\begin{theorem} (see, e.g., \cite{A1}) There exist $\e_0, c_1$ such that for $0<\e<\e_0$ we have 
$$| I(p(t),q(t),\e t)-I(p(0),q(0),0)| <c_1\e \;\; {\rm for} \;\;\;0\le t\le 1/\e\,.$$
\end{theorem}

Now assume that for each   $\tau={\rm const } \in [0,1]$, and each $I_*\in (a_1, b_1)$ Hamiltonian $H$ (\ref{ad_class1})  has a unique trajectory  with the action $I=I_*$.
Consider the corresponding  quantum system (\ref{S_h}).   The  operator $\cH_{\tau}$  has a series of eigenfunctions $\vp_{s}(\tau)=\vp_s(\tau,x)$ such that 
\begin{equation}
\label{norm}
||\vp_{s}(\tau)||=1, \qquad  \vp_{s}(\tau,x)\to0\;\;\text{as}\;\; x\to\infty,
\end{equation}
and the corresponding 
eigenvalues are
 $\lambda_s(\tau) = E(I_s, \tau)+ O(\hbar^2)$, where $I_s=\hbar(s+1/2) \in [a_1, b_1]$  (this is the 
    Bohr-Sommerfeld quantisation rule, see  \cite{MF}).  We assume that  $V$ is such that  
    the convergence to zero in \eqref{norm}  is 
    faster than  $|x|$ in any negative  power.     Let  $u(t,x)$ be a solution 
    of non-stationary equation   (\ref{S_h}) with a pure state initial condition  $u(0,x)=\vp_{s_0}(0)$.   
    Denote by     $\Prr_{(\alpha, \beta)}^\tau$ the orthogonal projector in   $L^2(\R)$ onto the  linear span
     of vectors $\vp_{s}(\tau)$  with $I_s\in (\alpha, \beta)$. The approach in \cite{Bor} leads to the following    
\begin{conjecture} 
\label{1d_quantum adiabatic} 
 There exist $ \e_0, c_1$ such that  if $0<\e<\e_0$ and  $0<\hbar\le \e$, then  for any $m\ge1$ and a suitable 
  $c_2(m)>0$ we have 
\begin{equation}
\label{qad_est}
\sup_{0\le t\le \e^{-1}   }||u(t)-\Prr_{(I_{s_0}-c_1\e, I_{s_0}+c_1\e)}^{\e t}  u (t)|| <c_2(m) \left(\frac{\hbar}{\e}\right)^m \,.
\end{equation}
\end{conjecture}

Thus $u(t,\cdot)$ stays close  to the eigenspace that  corresponds to eigenvalues from $O(\e)$-neighbourhood of $E( I_{s_0},\e t).$
\subsection{Systems with several  degrees  of freedom}
Now let  classical Hamiltonian  (\ref{ad_class1}) has $d>1$ degrees of freedom. As before, we
 assume that $V\in C^{\infty}$.  For each $\tau={\rm const}$ let the corresponding Hamiltonian system be completely integrable and in its phase space there is a domain filled by invariant tori. In this domain we  introduce action-angle variables $I=I(p,q,\tau)$,  $\ \chi=\chi(p,q,\tau) \in \T^d$. Invert these relations: $p=p(I, \chi, \tau), q=q(I, \chi, \tau)$. Suppose that there is a compact domain $ {\cal A}\Subset \R^d_+$ such that
  the map $(I, \chi, \tau) \mapsto (p, q, \tau)$ is smooth for  $I\in {\cal A},\chi\in \T^d, \tau\in[0,1] $.
 We  express Hamiltonian (\ref{ad_class1}) via the action variables and slow time, $H(p,q,\tau)=E(I, \tau)$, and
 denote by  $\omega(I,\tau)= \p E/\p I$ the frequency vector  of the unperturbed motion. We assume that the system is non-degenerate or iso-energetically nondegenerate (see definition in \cite{A1}, Appendix 8).   The dynamics of the variables 
  $(I, \chi)(t)=(I,\chi)(p(t), q(t), \e t)  $ is described by a Hamiltonian of the form (see  \cite{A1}, Sect. 52F)
\begin{equation}
\label{H_I}
{\cal H}(I,\chi,\tau,\e)=E(I, \tau)+\e H_1(I,\chi,\tau),
\end{equation}
where $H_1$ is  a smooth function on ${\cal A}\times\T^d\times[0,1]$.

Let $K_0$ be a compact set  in $ \R^{2d}$. For $(p_0,q_0)\in K_0$ denote by
$(p,q)(t)= (p,q)(t,p_0,q_0)$  a solution of the perturbed system 
with initial condition $(p,q)(0)=(p_0, q_0)$.
\begin{theorem}  (see, e.g., \cite{AKN, LochM}). If $0<\e<\e_0$, then
$$
\int\limits_{K_0}  \sup_{0\le t\le \e^{-1}   }| I(p(t),q(t),\e t)-I(p(0),q(0),0)|dp_0dq_0 <c_1\sqrt{\e}       \,.
$$
\end{theorem} 
In systems with $d>1$ degrees of freedom the value of  action-vector as a function of time 
may change considerably
 for some initial conditions due to the effect of resonance between   unperturbed frequencies, i.e. components of the vector $\omega(I,\tau)$. We  say that there is a resonance for some $(I, \tau)$ if $(k\cdot\omega)(I,\tau)=0$ for a suitable 
 vector   $k\in\Z^d\setminus\{0\}$ (here $\cdot$  denotes the  Euclidian scalar product). 

Now consider  corresponding  quantum system (\ref{S_h}). Under some conditions,  the  operator $\cH_{\tau}$ 
 has a series of eigenfunctions $\vp_{s}(\tau)=\vp_s(\tau,x)$, $s\in\Z^d$,  satisfying \eqref{norm},  with eigenvalues $\lambda_m(\tau) = E(I_m, \tau)+ O(\hbar^2)$, where $I_m=\hbar(m+\frac{1}{4}\kappa) \in {\cal A}$,  $m\in \Z^d_+$, and  $\kappa\in\Z^d$ is the vector of the  Maslov-Arnold indices \cite{MF}
  (the Bohr-Sommerfeld quantisation rule).  Consider now the solution $u(t,x)$ of non-stationary equation   (\ref{S_h}) with a pure state  initial condition  $u(0,x)=\vp_{m_0}(0)$.
  If we  fix some small  $\hbar $ and proceed to the limit as $\e\to 0$, then Theorem \ref{tB-F} would apply.  However, now we are interested  in another limit, when a small 
  $\e$ is fixed and $\hbar\to 0$.  Not much is known about the corresponding limiting dynamics. So we will formulate 
 natural  {\it hypotheses }  about  the limiting quantum dynamics as $\hbar\to0$ and  will use them jointly with the 
 known   results about dynamics for  classical Hamiltonian (\ref{ad_class1}) with small $\e$. 
  
For Theorem ~\ref{tB-F}  to hold 
it is important that  $\lambda_{m_0}(\tau)$ is an isolated eigenvalue for all $\tau$. Consider the
distance between $\lambda_{m}(\tau)$ and $\lambda_{m_0}(\tau)$, where 
$ m, m_0\in\Z^d$ are such that 
$m\ne m_0$ and  $|m-m_0|\sim 1$:  
\begin{equation*}\begin{split}
\lambda_{m}(\tau)-\lambda_{m_0}(\tau)&=  E(I_m, \tau)-E(I_{m_0}, \tau)  + O(\hbar^2)\\
&=  (I_m-I_{m_0})  \cdot \omega(I_{m_0},\tau) +O((I_m-I_{m_0})^2)+ O(\hbar^2)\\
&=\hbar(m-{m_0})\cdot\omega(I_{m_0},\tau)+ O(\hbar^2) .
\end{split}
\end{equation*}
Thus if there is no resonance at  $(I_{m_0}, \tau)$, then distance between $\lambda_{m_0}(\tau)$ and nearby eigenvalues is $\sim \hbar$. However, if there is a resonance $k\cdot\omega(I_{m_0},\tau)=0$, then  $\lambda_{m_0+\nu k}(\tau)-\lambda_{m_0}(\tau)=O(\hbar^2)$ for integer $\nu\sim1$. Thus classical resonances correspond to almost multiple points of the spectrum of the quantum problem.  Therefore  it seems that they  should 
also manifest themselves in the  quantum adiabaticity. 
\medskip

For Hamiltonian  (\ref{ad_class1})  there is a rather detailed information about dynamics in the two-frequency case
 $d=2$.   We will now use this information and the Bohr-Sommerfeld quantisation rule  to 
 state some conjectures   about dynamics for the  2d~quantum  system   (\ref{S_h}).

Following P. Dirac \cite{Dir25} we  assume that \footnote[4] {Condition (\ref{A}) just means that the ratio of frequencies changes with non-zero rate along solutions  of the system with Hamiltonian (\ref{H_I}): $\omega_2^2\frac{d}{dt}\big(
\frac{\omega_1}{\omega_2}
\big) >c^{-1}\e$. Similarly, condition (\ref{barA}) means that ratio of frequencies changes with non-zero rate in adiabatic dynamics: $\omega_2^2\frac{d}{dt}\big(
\frac{\omega_1}{\omega_2}
\big)_{I={\rm const}} >c^{-1}\e$. }
\begin{equation}
\label{A}
\omega_2\frac{\p\omega_1}{\p \tau}-\omega_1\frac{\p\omega_2}{\p \tau}-\left(\omega_2\frac{\p\omega_1}{\p I}-\omega_1\frac{\p\omega_2}{\p I}\right)\frac{\p H_1}{\p\chi}>c^{-1}
\end{equation}
for all $I,\vp$. General result by V.~I.~Arnold about averaging in two-frequency systems \cite{Arn65, AKN} implies  that in this case 
\begin{equation}
\label{Aest}
| I(p(t),q(t),\e t)-I(p(0),q(0),0)| <c_1\sqrt {\e} \quad {\rm for} \ 0\le t\le 1/\e\,.
\end{equation}
On the basis of the  Bohr-Sommerfeld quantisation rule and by analogy with Conjecture~\ref{1d_quantum adiabatic}
it is natural to conjecture that for $0\le t\le 1/\e$ 
the total probability $|u(t)|_{L_2}^2$ is mostly concentrated in the states, 
corresponding to actions from the $C\sqrt\e$-vicinity of the original action $I_{s_0}$.

Now assume that instead of (\ref{A}) the following condition is satisfied (cf. the forth footnote)
\begin{equation}
\label{barA}
\omega_2\frac{\p\omega_1}{\p \tau}-\omega_1\frac{\p\omega_2}{\p \tau}>c^{-1}\,.
\end{equation}
This  is a particular case of a condition introduced by V.~I.~Arnold in \cite{Arn65}.
If, in addition to (\ref{barA}), 
some general position condition is satisfied  (see details in \cite{AKN}), then   estimate (\ref {Aest}) in which 
 $\sqrt{\e}$ is replaced with $\sqrt{\e}|\ln\e|$ holds for all initial data outside 
 a set of measure
  $O(\sqrt{\e})$ \cite{AKN}, Sect. 6.1.8. The later  set mainly consists  of initial data for trajectories with {\it capture into resonance},
  along these  trajectories actions change by values $\sim 1$. Since for some initial 
  data  $I(0), \chi(0)$ the solution $I(t)$ is not localised in the vicinity of $I(0)$, then   we should 
    not expect for the  quantum system (\ref{S_h}) any estimate similar to that of Conjecture~\ref{1d_quantum adiabatic}, where the amplitudes of eigenmodes  tend to 0 as $\hbar \to 0$ outside some small interval of actions.  

Consider  classical Hamiltonian   (\ref{ad_class1})  under condition (\ref{barA}). Then the capture is only
possible for a finite number of resonances, 
and the dynamics with a capture into resonance   $k_1\omega_1+k_2\omega_2=0$ with 
co-prime $k_1, k_2$ is the following  \cite{N05}.  Denote $(I, \chi)(t)=(I,\chi)(p(t), q(t), \e t)  $.
  Suppose that at  the initial moment  $t=0$ we have no resonance:
 $$
 k_1\omega_1(I(0),0)+k_2\omega_2(I(0),0)\ne0,
 $$
 and let $\tau_*\in(0,1)$ be the first moment when the resonance occurs:

 $$
 k_1\omega_1(I(0),\tau_*)+k_2\omega_2(I(0),\tau_*)=0.
 $$ 
  Then for $0\le\e t\le \tau_*$ the values of actions are approximately
 conserved:
 $$
 I(t)=I(0)+O(\sqrt{\e}\ln\e) \,. 
 $$
 For $\tau_*\le\e t\le1$ the system is captured into resonance, and evolution of actions is described by two relations:
 \begin{equation*}\begin{split}
&  k_1\omega_1(I(t),\e t)+k_2\omega_2(I(t),\e t )=O(\sqrt{\e}\ln\e),\\
&   k_2 I_1(t)- k_1 I_2(t)=  k_2 I_1{(0)}- k_1 I_2{(0)}+ O(\sqrt{\e}\ln\e)\;.\;  
  \end{split}
\end{equation*}
First of them means that the system stays near the resonance, while the second  says that the dynamics has an approximate  first  integral. Jointly the two relations approximately 
define the trajectory $I(t)$ for $\tau_*\le \e t\le1$.

  Based on this description and the Bohr-Sommerfeld quantisation rule, by analogy with 
  Conjecture~\ref{1d_quantum adiabatic} we conjecture that for the quantum problem  \eqref{S_h} the  
  capture in resonance in the classical system  \eqref{ad_class1} results in transfer of an $C\e$-amount of the total probability from the vicinity of the   initially excited pure state, corresponding to the action $I_{s_0}$, to the
   vicinity of a state $s_t\in \Z^2$ such that the lattice vector
   $I(t)=\hbar(s_t+\frac{1}{4}\kappa) $   satisfies the two relations above. This transfer happens for  $t\ge\e^{-1}\tau_*$. When 
   $\hbar\to0$, this  $C\e$-amount stays positive of order $\e$.
   \smallskip

For the dynamics of phases of 
captured into resonances points also there is a  more detailed description
 \cite{N05}. Consider
 the resonant phase $\gamma= k_1\chi_1+k_2\chi_2$. It turns out that  the  behaviour of  $\gamma$ is described by an auxiliary Hamiltonian system with one degree of freedom and the  Hamiltonian of the form
$$
F=\sqrt{\e}\left(\alpha(\tau)p_{\gamma}^2/2 +f(\gamma, \tau) +L(\tau)\gamma \right)\,.
$$ 
Here $p_{\gamma}, \gamma$ are canonically conjugate variables, function $f$ is $2\pi$-periodic in $\gamma$, 
and $\alpha,L\ne 0$.
In the phase portrait of the system for frozen $\tau$ there are domains of oscillations of $\gamma$. Motion in these 
domains can be approximately represented as composition of motion along a  trajectory  of  Hamiltonian $F$ with  frozen $\tau$  and slow evolution of this trajectory  due to the  change of $\tau$. 
This evolution follows the  adiabatic rule: the area surrounded by the  trajectory remains  constant. 
In the original variables $p,q$ this motion is represented as a motion along slowly evolving torus. 
Angular variables on this torus are $\gamma$ and $\psi=l_1\vp_1+l_2\vp_2$, where $l_1$ and $l_2$ are integers 
 such that $k_1l_2-k_2l_1=1$.  This torus drifts along the resonant  surface $k_1\omega_1+k_2\omega_2=0$ as it was described above. It is not known  which   quantum object corresponds to it. 

\bibliography{kam_q}
\bibliographystyle{amsalpha}
\end{document}